\documentclass[11pt]{article}
\usepackage{amsfonts}

\usepackage{graphics}

\usepackage{indentfirst}
\usepackage{cite}
\usepackage{latexsym}
\usepackage{amsmath}
\usepackage{amssymb}
\usepackage[dvips]{epsfig}
\usepackage{amscd}

\usepackage{amssymb,mathrsfs,psfrag,graphicx}
\usepackage{hyperref}

\newtheorem{theorem}{Theorem}[section]
\newtheorem{remark}{Remark}[section]

\newtheorem{definition}{Definition}[section]
\newtheorem{lemma}[theorem]{Lemma}

\newcommand{\n}{\rho}

\newcommand{\mr}{\mathbb{R}}
\newcommand{\lm}{\lambda}

\renewcommand{\div}{ {\rm div }  }

\newcommand{\pa}{\partial}
\renewcommand{\r}{\mathbb{R}}

\renewcommand{\b}{B_{N_*}}

\newcommand{\ia}{\int_0^T}

\newcommand{\bt}{\begin{theorem}}
\newcommand{\bl}{\begin{lemma}}
\newcommand{\el}{\end{lemma}}
\newcommand{\et}{\end{theorem}}
\newcommand{\ga}{\gamma}

\newcommand{\al}{\alpha}
\newcommand{\de}{\delta}
\newcommand{\ve}{\varepsilon}
\newcommand{\la}{\label}

\newcommand{\ol}{\overline}

\newcommand{\bn}{\begin{eqnarray}}
\newcommand{\en}{\end{eqnarray}}
\newcommand{\bnn}{\begin{eqnarray*}}
\newcommand{\enn}{\end{eqnarray*}}

\newcommand{\bnnn}{\begin{eqnarray*}}
\newcommand{\ennn}{\end{eqnarray*}}

\newcommand{\ba}{\begin{aligned}}
\newcommand{\ea}{\end{aligned}}
\newcommand{\be}{\begin{equation}}
\newcommand{\ee}{\end{equation}}
\def\O{{\r^2 }}
\def\p{\partial}
\def\norm[#1]#2{\|#2\|_{#1}}

\def\o{\omega}

\def\la{\label}

\def\na{\nabla}

\makeatletter      
\@addtoreset{equation}{section}
\makeatother

\makeatletter      
\@addtoreset{equation}{section}
\makeatother       

\title{A blow-up criterion of strong solutions to the 2D
compressible magnetohydrodynamic equations }

\author{ Teng W{\small ANG}\thanks{Institute of Applied Mathematics, AMSS, Chinese Academy of Sciences, Beijing 100190, China
({\tt tengwang@amss.ac.cn})}
 }

\date{ }

\begin{document}
\maketitle

\begin{abstract} This paper establishes a blow-up criterion of strong solutions to
the two-dimensional compressible magnetohydrodynamic (MHD) flows.
The criterion depends on the density, but is independent of the velocity and the magnetic field.
More precisely,  once the strong solutions blow up, the $L^{\infty}$-norm for the density tends to infinity.
In particular, the vacuum in the solutions is allowed.

\end{abstract}


\section{Introduction}

We consider the system  for the two-dimensional viscous,
compressible magnetohydrodynamic (MHD) flows in the Eulerian coordinates as follows,
\be \label{a1}  \begin{cases} \n_t + \div(\n u) = 0,\\
 (\n u)_t + \div(\n u\otimes u) + \nabla P(\n) = \mu\Delta u + (\mu + \lambda)\nabla(\div u)+(\na\times H)\times H,\\
 H_t-\na \times(u\times H)=\nu\Delta H,\quad \div H=0,
\end{cases}\ee
where $\n=\n(x,t),$ $u=(u^1, u^2)(x,t),  $ $H=(H^1,H^2)(x,t)$,  and \be P(\n)=R\n^\ga\,\,( R>0, \ga>1)  \ee
are the fluid density, velocity, magnetic field and pressure, respectively. $R$ is a positive constant.
Without loss of generality, we assume that $R=1$. The constant
viscosity coefficients $\mu$ and $\lambda$ satisfy the physical
restrictions: \be\la{h3} \mu>0,\quad \mu +\lambda\ge 0.
\ee
The constant $\nu>0$ is the resistivity coefficient which is inversely proportional to the
electrical conductivity constant and acts as the magnetic diffusivity of magnetic fields.
We consider  the Cauchy problem of (\ref{a1}) with the initial data
\be\label{h2}
(\n,u,H)|_{t=0}:=(\n_0,u_0,H_0), \quad x\in  \r^2,
\ee
and the boundary condition at the far fields
\be\label{far}
(\n,u,H)\rightarrow(0,0,0), \quad {\rm as}~~ |x|\rightarrow\infty.
\ee

There have been large literature on the compressible MHD system \eqref{a1} by many physicists and mathematicians
due to its physical importance, complexity, rich phenomena and mathematical challenges, see
\cite{cw,cw2,fjn,fy2,he-xin,hw1,hw2,ko,ka,lxz,chl,lv-shi-xu,uks,w,xu-zhang} and the references therein.
When the initial density is uniformly positive, the local existence of strong solutions to the three-dimensional
compressible MHD was proved by Vol'pert-Khudiaev \cite{vk}. Then,  Kawashima \cite{ka}
firstly obtained the global existence when the initial data are close to a non-vacuum equilibrium in $H^3$-norm.
In particular, the theory requires that the solutions have small oscillations around a uniform
non-vacuum state so that the density is strictly away from the vacuum.

For general large initial data, the global well-posedness of classical solutions to the compressible MHD system remains
completely open. One of the main difficulties is that even the initial density is absence of vacuum, one could not know whether the vacuum
states may occur or not within finite time. The system \eqref{a1} may degenerate in the presence of vacuum, which
produces new difficulty in mathematical analysis. Therefore, it is interesting  to study the
solutions with vacuum for the compressible MHD system.
In fact, the local well-posedness of strong solutions of the system \eqref{a1} in three dimensions was established by
Fan-Yu \cite{fy2} when the initial density may contain vacuum. Recently, Li-Xu-Zhang \cite{lxz}
proved the global existence of classical solutions with vacuum as far field condition to \eqref{a1}-\eqref{far}
in $\mathbb{R}^3$ for regular initial data with small energy but possibly large oscillations. However, the two
dimensional case is quiet different from the three dimensional case when the far field condition is vacuum.
Precisely speaking, the difference between 2D and 3D is that, if a function $u$ satisfies $\na u\in L^2(\mathbb{R}^2)$, it is impossible to imply that $u\in L^p(\mathbb{R}^2)$,
for any $p>1$, while if $\na u\in L^2(\mathbb{R}^3)$, then $u\in L^6(\mathbb{R}^3)$. More recently,
Lv-Huang \cite{chl} proved the local existence of strong solutions for the Cauchy problem of the two-dimensional
compressible MHD equations with vacuum as far field density by weighted energy estimate and
Lv et al. \cite{lv-shi-xu} generalized the previous work of Li-Xu-Zhang \cite{lxz} on two-dimensional case.

Until now, all the global existence of solutions for the compressible MHD equations were obtained with some "smallness"
assumptions on the initial data. However, there still remains a longstanding open problem: whether the strong (or classical) solutions to
the compressible MHD system \eqref{a1} can exist globally or not?
Thus, it is important to study the mechanism of blow up and structure of possible singularities of
strong solutions to the compressible MHD system \eqref{a1}. Along this direction, Xu-Zhang in \cite{xu-zhang} established the
following Serrin's type criterion to the three-dimensional isentropic compressible MHD equations:
\be\label{l-serrin-2}
\lim_{T\rightarrow T^*}\Big(\|\n\|_{L^{\infty}(0,T;L^{\infty})}
+\|u\|_{L^s(0,T;L^r_w)}\Big)=\infty,
\ee
where $T^*\in(0,\infty)$ is the maximal time of existence for
strong (or classical) solutions $(\n, u)$, $L^r_w$ denotes the weak $L^r$-space
and $r$, $s$ satisfies
\be\label{serrin}
\frac{2}{s}+\frac{3}{r}\leq 1, \quad 3<r\leq \infty.
\ee
Then Huang-Li in \cite{hl} extended  the result of \cite{xu-zhang} to the non-isentropic case.


Although the Serrin type criterion for the Cauchy problem of three-dimensional MHD flows
has been well established by Xu-Zhang in \cite{xu-zhang}. However, the two-dimensional case
with the initial density containing vacuum becomes more difficult since the analysis of \cite{xu-zhang} for 3D case
depends crucially on the $L^6$-bound on the velocity, while in 2D case, the velocity may not belong to $L^p(\mathbb{R}^2)$
for any $p>1$.  The main aim of this paper is to establish a blow-up criterion for the two-dimensional compressible MHD system.

Before stating the main results, we first explain the notations and
conventions used throughout this paper. For $R>0$, set
$$B_R  \triangleq\left.\left\{x\in\r^2\right|
\,|x|<R \right\} , \quad \int fdx\triangleq\int_{\r^2}fdx.$$ Moreover,
for $1\le r\le \infty$ and $k\ge 1$, the standard homogeneous
and inhomogeneous Sobolev spaces are defined as follows:
   \bnn  \begin{cases}L^r=L^r(\r^2 ),\quad D^{k,r}=D^{k,r}(\r^2)
   =\{v\in L^1_{\rm loc}(\r^2)| \na^k v\in L^r(\r^2)\}, \\ D^1=D^{1,2},\quad
W^{k,r}  = W^{k,r}(\r^2) , \quad H^k = W^{k,2} . \\
 \end{cases}\enn

 Next, we give the definition of strong solution to \eqref{a1} as follows:
\begin{definition}[strong solution]\label{def} $(\rho, u,H)$ is called a
strong solution to \eqref{a1} in $\mathbb{R}^2\times(0,T)$,
if for some $q_0>2$ and $a>1$,
\be \label{strong}  \begin{cases}
\n\geq 0,\quad \n\bar x^{a}\in C([0,T];L^1\cap H^1\cap W^{1,q_0}),
\quad\n_t\in C([0,T];L^2\cap L^{q_0}),\\
 (u,H)\in C([0,T];D^1\cap D^2)\cap L^2(0,T;D^{2,q_0}),
 \quad H\in C([0,T];H^2),\\
 (\n^{1/2}u_t,H_t)\in L^{\infty}(0,T;L^2),\quad
 (u_t,H_t)\in L^2(0,T;D^1),
\end{cases}\ee
and $(\n,u,H)$ satisfies \eqref{a1} $a.e.$ in $\mathbb{R}^2\times(0,T)$,
where
  \be\label{2.07} \bar x\triangleq(e+|x|^2)^{1/2} \log^2 (e+|x|^2) .\ee

\end{definition}

Without loss of generality, assume that the initial density $\n_0$ satisfies
\be\la{oy3.7} \int_{\r^2} \n_0dx=1,\ee  which implies that there exists a
positive constant $N_0$ such that
\be\label{oy3.8} \int_{B_{N_0}}  \n_0  dx\ge \frac12\int\n_0dx=\frac12.\ee

Our main result can be stated as follows:
\begin{theorem}\label{th1} Let $\Omega=\mathbb{R}^2$. In addition to \eqref{oy3.7} and  \eqref{oy3.8},
suppose that  the initial data $(\n_0,u_0,H_0)$ satisfy, for any given numbers $a>1$, $q>2$,
\begin{equation} \label{co1}
\left\{\begin{array}{lll}
\n_0\geq 0, \quad  \bar x^{a} \rho_0\in   L^1\cap H^1\cap W^{1,q},
\quad \n_0u_0^2+\n_0^{\gamma} \in L^1,\\
 (u_0,H_0)\in D^1\cap D^2,\quad \bar{x}^{a/2}H_0\in L^2,
\quad \bar{x}^{a/2}\na H_0\in L^2,
\end{array}\right.
\end{equation}
 and the compatibility condition
 \be\label{cc}
 -\mu\Delta u_0-(\mu+\lambda)\nabla\div u_0 +\nabla P(\n_0)-(\nabla\times H_0)\times H_0
 =\n_0^{1/2}g,
 \ee holds for some $g\in L^2$.

  Let $(\n,u,H)$ be a strong solution to the Cauchy problem \eqref{a1}, \eqref{far} and \eqref{h2},
  satisfying \eqref{strong}. If $T^*<\infty$ is the maximal time of existence, then
  \be\label{blowup}
  \lim_{T\rightarrow T^*}\|\n\|_{L^{\infty}(0,T;L^{\infty})}=\infty.
  \ee

\end{theorem}


\begin{remark} \label{re2}
Theorem \ref{th1} means the blow-up criterion \eqref{blowup} is independent of the velocity $u$ and the magnetic field $H$.
\end{remark}

We now make some comments on the analysis of this paper. As previously mentioned, the $L^6$-norm of velocity $u$
plays crucial role in 3D MHD system, while for 2D, $u$ may not belong to $L^p(\mathbb{R}^2)$, for any $p>1$. The key observation
of this paper is that, if we restriction the initial data in a smaller space, i.e. $\bar{x}^a\rho_0\in L^1(\mathbb{R}^2)$
for some positive constant $a>1$ (see \eqref{co1}),
we can show that $u\bar x^{-\eta}$ belongs to $L^{p_0}(\mathbb{R}^2)$, for some positive constant $p_0>1$,
and $\eta\in(0,1]$ (see \eqref{u-wei'}), here $\bar{x}=(e+|x|^2)^{1/2}\log^2(e+|x|^2)$.
For this , we need to manipulate the weighted energy estimates throughout the Section \ref{se3} below.

Before finishing the introduction, we recall some related works. The global existence of weak solutions to the system \eqref{a1}
 was established by Hu-Wang \cite{hw1,hw2}.
 If there is no electromagnetic field effect, \eqref{a1} turns to be the  compressible Navier-Stokes equations.
 For the blow-up criterion of the compressible Navier-Stokes equations,
 we refer to \cite{hx2,hlx,hlw,sun-zhang,swz-1} and the references therein.

The rest of the paper is organized as follows: In the next section, we
collect some elementary facts and inequalities for the blow-up analysis.
The main result, Theorem \ref{th1}, is proved in Section \ref{se3}.

\section{Preliminaries}\label{se2}

In this section, we will recall some  known facts and elementary
inequalities which will be used frequently later.

The local existence of strong solutions when the initial density may not be
positive can be proved in a similar way as in \cite{li-liang} (cf. \cite{chl}).

\begin{lemma}   \label{th0}  Assume  that the initial data
 $(\n_0,u_0,H_0)$ satisfy \eqref{co1}. Then there exist  a small time
$T_1 >0$  and a unique strong solution $(\rho , u,H )$ in the sense of
Definition \ref{def} to the Cauchy problem   \eqref{a1}-\eqref{h2} in
$\O\times(0,T )$.
 \end{lemma}

Next, the following well-known Gagliardo-Nirenberg inequality (see \cite{nir})
will be used later.

\begin{lemma}
[Gagliardo-Nirenberg]\label{l1} For  $p\in [2,\infty),q\in(1,\infty), $ and
$ r\in  (2,\infty),$ there exists some generic
 constant
$C>0$ which may depend  on $p,q, $ and $r$ such that for   $f\in H^1({\O }) $
and $g\in  L^q(\O )\cap D^{1,r}(\O), $    we have \be
\la{g1}\|f\|_{L^p(\O)}^p\le C \|f\|_{L^2(\O)}^{2}\|\na
f\|_{L^2(\O)}^{p-2} ,\ee  \be
\la{g2}\|g\|_{C\left(\ol{\O }\right)} \le C
\|g\|_{L^q(\O)}^{q(r-2)/(2r+q(r-2))}\|\na g\|_{L^r(\O)}^{2r/(2r+q(r-2))} .
\ee
\end{lemma}

The following weighted $L^p$ bounds for elements of the Hilbert space $  D^{1}(\O)  $ can be found in \cite[Theorem B.1]{L2}.
\begin{lemma} \label{1leo}
   For   $m\in [2,\infty)$ and $\theta\in (1+m/2,\infty),$ there exists a positive constant $C$ such that we have for all $v\in  D^{1,2}(\O),$ \be\la{3h} \left(\int_{\O} \frac{|v|^m}{e+|x|^2}(\log (e+|x|^2))^{-\theta}dx  \right)^{1/m}\le C\|v\|_{L^2(B_1)}+C\|\na v\|_{L^2(\O) }.\ee
\end{lemma}

The following lemma was deduced in \cite{lx1}, we only state it here without proof.

\begin{lemma}\label{lemma2.6} For $\bar x$   as in \eqref{2.07},
suppose that $\n  \in L^\infty(\O)$ is a   function such that
\be \la{2.12}  0\le \n\le M_1, \quad M_2\le \int_{\b}\n dx ,\quad \n \bar x^\alpha \in L^1(\r^2),\ee
for $ N_*\ge 1 $ and positive constants $   M_1,M_2, $  and   $\al.$  Then, for $r\in [2,\infty),$ there exists a positive constant $C$ depending only on $  M_1, M_2, \alpha,   $ and $ r$  such that
 \be\label{z.1}\left(\int_{\r^2}\n |v |^r dx\right)^{1/r}  \le C  N_*^3  (1+\|\n\bar x^\al\|_{L^1(\r^2)})  \left(  \|\n^{1/2} v\|_{L^2(\r^2)} + \|\na  v \|_{L^2(\r^2)}\right) ,\ee for each $v\in \left.\left\{v\in D^1 (\O)\right|\n^{1/2}v\in L^2(\r^2) \right\}.$

\end{lemma}

Next, for $ \nabla^{\perp}\triangleq (-\p_2,\p_1)$,
denoting the material derivative of $f $ by $\dot f\triangleq f_t+u\cdot\nabla f$.
We state some elementary estimates which follow from
(\ref{g1}) and the standard $L^p$-estimate  for the following elliptic
system derived from the momentum equations in \eqref{a1}$_2$:
\be\label{F}
\triangle F = \text{div}\left(\n\dot{u}-\text{div}(H\otimes H-\frac12|H|^2\mathbb{I}_2)\right),
\ee
and
\be\label{w}
\mu \triangle \o =
\nabla^\perp\cdot\left(\n\dot{u}-\text{div}(H\otimes H-\frac12|H|^2\mathbb{I}_2)\right) ,
\ee
where
\be\label{evf}
F\triangleq(2\mu+\lambda)\div u-P(\n),\quad \o=\partial_1u^2-\partial_2u^1.
\ee
The symbol $\otimes$ denotes the Kronecker tensor product, i.e.
$H\otimes H=(H^iH^j)_{2\times 2}$ and "$\mathbb{I}_2$" denotes the
$2\times 2$ unit matrix.

\begin{lemma} \label{le2.5}
  Let $(\n,u,H)$ be a smooth solution of
   (\ref{a1}).
    Then for   $p\ge 2$ there exists a   positive
   constant $C$ depending only on $p,\mu,$ and $\lambda$ such that
\begin{eqnarray}
    &&\|{\nabla F}\|_{L^p(\O)} + \|{\nabla \o}\|_{L^p(\O)}
   \le C(\norm[L^p(\O)]{\n\dot{u}}+\norm[L^p(\O)]{|H||\na H|}),\label{h19}\\
&&\norm[L^p(\O)]{F} + \norm[L^p(\O)]{\o}
   \le C \left(\norm[L^2(\O)]{\n\dot{u}}+\norm[L^2(\O)]{|H||\na H|}\right)^{1-2/p}\nonumber\\
   &&\quad\cdot\left(\norm[L^2(\O)]{\nabla u} + \norm[L^2(\O)]{P }\right)^{2/p} ,
\label{h20}
\\
  &&\norm[L^p(\O)]{\nabla u} \le C \left(\norm[L^2(\O)]{\n\dot{u}}+\norm[L^2(\O)]{|H||\na H|}\right)^{1-2/p }\nonumber
   \\&&\quad\cdot\left(\norm[L^2(\O)]{\nabla u}
   + \norm[L^2(\O)]{P }\right)^{2/p}+ C\norm[L^p(\O)]{P}. \label{h18}
\end{eqnarray}
\end{lemma}
{\bf Proof}: On the one hand, by the standard $L^p$-estimate for the elliptic
systems (see \cite{nir}), (\ref{F}) \eqref{w} yield (\ref{h19}) directly, which, together
with (\ref{g1}), \eqref{F} and \eqref{w}, gives (\ref{h20}).
On the other hand, since $-\Delta u=-\na {\rm div}u -\na^\perp\o,$ we have \bn\la{kq1}\na u=-\na(-\Delta)^{-1}\na {\rm
div}u-\na(-\Delta)^{-1}\na^\perp \o.\en Thus applying the standard $L^p$-estimate to \eqref{kq1} shows  \bnn \ba \|\na u\|_{L^p(\O)}&\le C(p) (\|{\rm
div}u\|_{L^p(\O)}+\|\o\|_{L^p(\O)})\\ &\le C (p)  \norm[L^p(\O)]{F} +C(p)  \norm[L^p(\O)]{\o} +
   C(p)  \norm[L^p(\O)]{P },\ea  \enn which,
along with (\ref{h20}), gives (\ref{h18}). Then, the proof of Lemma \ref{le2.5} is completed.

Finally,    the following Beale-Kato-Majda type inequality,
which was proved in \cite{hlx}, will be
used later to estimate $\|\nabla u\|_{L^\infty}$ and
$\|\nabla\n\|_{L^2\cap L^q} (q>2)$.
\begin{lemma}   \la{le9}  For $2<q<\infty,$ there is a
constant  $C(q)$ such that  the following estimate holds for all
$\na u\in L^2(\O)\cap D^{1,q}({\O }),$
\bnn \label{ww7}\ba
\|\na u\|_{L^\infty({\O })}&\le C\left(\|{\rm div}u\|_{L^\infty({\O })}
+\|\o\|_{L^\infty({\O })} \right)\log(e+\|\na^2 u\|_{L^q({\O })})\\
&\quad+C\|\na u\|_{L^2(\O)}+C.
\ea\enn
\end{lemma}

\section{\label{se3} Proof of Theorem \ref{th1}}

Let $(\n,u,H)$ be a strong solution to the problem \eqref{a1}-\eqref{h2} as describe in Theorem \ref{th1}.
Suppose that \eqref{blowup} were false, that is,
\be\label{non-blowup}
  \lim_{T\rightarrow T^*}\|\n\|_{L^{\infty}(0,T;L^{\infty})}=M_0<\infty.
  \ee

First, the standard energy estimate yields
 \be \label{basic}\sup_{0\le t\le T}\int\left(\frac{1}{2}\n|u|^2
 +\frac{P}{\ga-1} +\frac{|H|^2}{2}\right)dx+\ia\int\left( \mu|\na
u|^2 +\nu|\na H|^2\right) dxdt\le C,
 \ee
for $0\leq T\leq T^*$. Throughout this paper, several positive generic constants
are denoted by $C$ and $C_i (i=1,2)$ depending only on
$M_0,~\mu,~\lambda,~\nu,~T^*,~q,~a$ and the initial data.

The following lemma is based on \eqref{basic}.
\begin{lemma}\label{le-he-xin}
It holds that for $q\in[2,\infty)$ and $0\leq T\leq T^*$,
\be\label{h1}
\|H\|_{L^{\infty}(0,T;L^q)}+\int_0^T\int|H|^{q-2}|\nabla H|^2 dxdt\leq C.
\ee
\end{lemma}

{\bf Proof}: We prove \eqref{h1} as in He-Xin \cite{he-xin}. Multiplying $\eqref{a1}_3$ by
$q|H|^{q-2}H$ and integrating the resulting equation over $\O$ yield that
\be\label{h1-1}
\ba
&\frac{\rm d}{\rm dt}\int|H|^qdx+\nu\int\left(q|H|^{q-2}|\nabla H|^2+q(q-2)|H|^{q-2}|\nabla|H||^2\right)dx\\
&\quad =q\int H\cdot\nabla u\cdot H|H|^{q-2} dx-(q-1)\int \div u |H|^q dx\\
&\quad \leq C\|\nabla u\|_{L^2}\||H|^{q/2}\|^2_{L^4}
       \leq C\|\nabla u\|_{L^2}\||H|^{q/2}\|_{L^2}\|\nabla|H|^{q/2}\|_{L^2}\\[2mm]
&\quad \leq \delta  \|\nabla|H|^{q/2}\|^2_{L^2}+C(\delta)\|\nabla u\|^2_{L^2}\||H\|^q_{L^q}.
\ea
\ee
Choosing $\delta$ suitable small in \eqref{h1-1}, we obtain \eqref{h1} directly after using Gronwall's
inequality and \eqref{basic}. Thus the proof of Lemma \ref{le-he-xin} is
completed.
\hfill $\Box$
\vspace{2mm}

Next, we give the key estimate on $\nabla u$ and $\nabla H$ in the following lemma.
\begin{lemma}\label{le1}
Under the condition \eqref{non-blowup}, it holds that for $0\leq T\leq T^*$,
\be\label{es-1}
\sup_{0\leq t\leq T}\Big(\|\nabla u\|^2_{L^2}+\|\nabla H\|^2_{L^2}\Big)
+\int_0^T\Big( \|\n^{1/2}\dot u\|^2_{L^2}+\|\nabla^2H\|^2_{L^2}\Big) dt\leq C.
\ee
\be\label{es-1'}
\int_0^T\Big(\|\nabla u\|^4_{L^4}+\|\nabla H\|^4_{L^4}\Big)dt\leq C.
\ee

\end{lemma}

{\bf Proof}: First, multiplying the momentum equation $\eqref{a1}_2$ by $\dot u$ and integrating the
resulting equation over $\O$ gives
\be\label{1-1}
\ba    \int  \n|\dot{u} |^2dx   &
= - \int  \dot{u}  \cdot\nabla Pdx + \mu \int \dot{u}\cdot\triangle
 u   dx + (\mu+\lambda)
 \int \dot{u}\cdot\nabla\text{div}u  \ dx \\
 &\quad-\frac{1}{2}\int\dot u\cdot\na |H|^2 dx
 +\int \dot u\cdot H\cdot\na H dx\triangleq\sum_{i=1}^5I_i,
\ea
\ee
where we have used the fact that
$$
(\nabla\times H)\times H=\div(H\otimes H)-\frac12\nabla|H|^2=H\cdot\nabla H-\frac12\nabla|H|^2.
$$
Since $P$ satisfies
\be\label{p}
P_t+\div(uP)+(\gamma-1) P\div u=0,
\ee
integrating by parts yields that
\be\label{i1}
\ba
-\int\dot u\cdot\nabla P dx &=\int\Big((\div u)_tP
-(u\cdot\nabla u)\cdot\nabla P\Big) dx\\
&=\left(\int \div u P dx\right)_t+\int\Big((\gamma-1)P(\div u)^2
+P\partial_i u\cdot\nabla u_i\Big) dx\\
&\leq \left(\int \div u P dx\right)_t+C\|\nabla u\|^2_{L^2}.
\ea
\ee
Integrating by parts also leads to
\be\label{i2}
\ba
\mu\int \dot{u}\cdot\Delta u dx&=-\frac{\mu}{2}\left(\|\nabla u\|^2_{L^2}\right)_t
-\mu\int\partial_i u^j\partial_i(u\cdot\nabla u^j)dx\\
&\leq -\frac{\mu}{2}\left(\|\nabla u\|^2_{L^2}\right)_t
+C\|\nabla u\|^3_{L^3}
\ea
\ee
and that
\be\label{i3}
\ba
(\mu+\lambda)\int\dot{u}\cdot\nabla\div u\ dx&=
-\frac{\mu+\lambda}{2}\left(\|\div u\|^2_{L^2}\right)_t
-(\mu+\lambda)\int \div u~\div(u\cdot\nabla u) dx\\
&\leq -\frac{\mu+\lambda}{2}\left(\|\div u\|^2_{L^2}\right)_t
+C\|\nabla u\|^3_{L^3}.
\ea
\ee
Using \eqref{a1}$_3$ and \eqref{basic}, we get
\be\label{i4}
\ba
I_4&=\frac{1}{2}\int |H|^2  \div u_tdx+\frac{1}{2}\int |H|^2  \div (u\cdot \na u)dx\\
&=\left(\int \frac{|H|^2}{2} \div udx\right)_t+\frac{1}{2}\int u\cdot\na |H|^2 \div u dx
+\frac{1}{2}\int |H|^2 \div(u\cdot \na u)dx\\
&\qquad-\int(H\cdot \na u+\nu \Delta H-H\div u)\cdot H\div udx\\
&=\left(\int \frac{|H|^2}{2} \div udx\right)_t-\frac{1}{2}\int |H|^2 (\div u)^2 dx+
\frac{1}{2}\int |H|^2 \pa_i u\cdot \na u^i dx\\
&\qquad-\int(H\cdot \na u+\nu \Delta H-H\div u)\cdot H\div u dx\\
&\leq \left(\int \frac{|H|^2}{2} \div udx\right)_t
+C(\ve)\int |H|^2|\na u|^2dx+\ve \|\na^2 H\|_{L^2}^2\\
&\leq \left(\int \frac{|H|^2}{2} \div udx\right)_t
+C(\ve)\|\na u\|_{L^3}^2\|H\|_{L^2}^{4/3}\|\na^2 H\|_{L^2}^{2/3}
+\ve \|\na^2 H\|_{L^2}^2\\
&\leq \left(\int \frac{|H|^2}{2} \div udx\right)_t
+C(\ve)\|\na u\|_{L^3}^3 +2\ve \|\na^2 H\|_{L^2}^2.
\ea\ee
Similarly, we have
\be\label{i5}
\ba
I_5\le - \left(\int H\cdot\na u\cdot H dx\right)_t
+ C(\ve)\|\na u\|_{L^3}^3  +2\varepsilon\|\na^2 H\|_{L^2}^2.
\ea
\ee

Putting \eqref{i1}-\eqref{i5} into  \eqref{1-1},
and recalling \eqref{non-blowup}, \eqref{h18} yields
\be\la{1-2}
\ba
 &B'(t)+\int  \n|\dot{u} |^2 dx\\
 &\le C \|\na u \|^2_{L^2}
 + C  (\ve)\|\nabla u\|_{L^3}^3 +4\varepsilon\|\na^2 H\|_{L^2}^2\\
&\le C(\ve)\left(\|\na u\|^4_{L^2}+1\right)
+C\||H||\na H|\|^2_{L^2}+\ve\|\n^{1/2}\dot u\|^2_{L^2}
+4\varepsilon\|\na^2 H\|_{L^2}^2.
\ea\ee
where
\be\ba \label{B(t)}
B(t)\triangleq &\frac{\mu  }{2}\|\nabla u\|_{L^2}^2
+\frac{ \lambda+\mu}{2}\|\text{div}u\|_{L^2}^2-\int  \text{div}u P dx\\
&-\frac{1}{2}\int\text{div}u|H|^2 dx+\int H\cdot\na u\cdot H dx.
\ea
\ee

Next, multiplying \eqref{a1}$_3$ by $ \triangle H$, and integrating by parts over $\mr^2$,  we have
\be\la{mm1}\ba
&\quad\frac{\rm d}{{\rm d}t}\int|\na H|^2 dx+2\nu\int|\nabla^2H|^2dx\\
&\quad\leq C\int |\na u||\na H|^2dx+ C\int |\na u|| H| |\Delta H|dx \\
&\quad\leq C\|\na u\|_{L^2} \|\na H\|^2_{L^4}
+C\|\na u\|_{L^2}\| H\|_{L^{\infty}}\|\na^2H\|_{L^2}\\
&\quad \le C\|\na u\|_{L^2} \|\na H\|_{L^2}\|\na^2 H\|_{L^2}
+C\|\na u\|_{L^2}\|\na^2 H\|^{3/2}_{L^{2}}\\
&\quad\leq \frac{\nu}{2}\|\na^2H\|^2_{L^2}+C\|\na u\|^4_{L^2}
+C\|\na u\|^2_{L^2}\|\na H\|^2_{L^2}.
\ea\ee

Choosing $C_1$ suitably large such that
\be \la{n2'}\ba&\frac{\mu }{4}\|\nabla u\|_{L^2}^2+ \|\na H\|_{L^2}^2 -C \|P\|_{L^2}^2\\
  & \le B(t)+C_1\|\na H\|_{L^2}^2\le  C \|\nabla u\|_{L^2}^2+C\|\na H\|_{L^2}^2+ C \|P\|_{L^2}^2,\ea\ee
 adding  \eqref{mm1} multiplied by $C_1$  to \eqref{1-2}, and choosing $\ve$ suitably small lead to
\be\label{n1}
\ba
&( B(t)+C_1\|\na H\|_{L^2}^2)'
+\frac12\int \left(\n|\dot{u} |^2 +\nu C_1|\nabla^2 H|^2\right)dx \\
& \le C+C\||H||\nabla H|\|^2_{L^2}+C(\|\na u\|^4_{L^2}+\|\na H\|^4_{L^2}).
\ea
\ee
Integrating \eqref{n1} over $(0,T)$, choosing $q=4$ in \eqref{h1},
and using \eqref{basic} and Gronwall's inequality, we obtain \eqref{es-1}.
We can get \eqref{es-1'} immediately by \eqref{h18} and \eqref{es-1}.

\hfill $\Box$
\vspace{2mm}


Next, we get some basic energy estimates on the magnetic field $H$.
\begin{lemma}\label{le2}
Under the condition \eqref{non-blowup}, it holds that for $0\leq T\leq T^*$,
\be\label{es-2}
\sup_{0\le t\leq T}\||H||\nabla H|\|^2_{L^2}+
\int_0^T\Big(\|\Delta |H|^2\|^2_{L^2}+\||\Delta H||H|\|^2_{L^2}\Big) dt
\leq C.
\ee
\end{lemma}

{\bf Proof}: We will follow an idea in \cite{lv-shi-xu}.
For $ a_1,a_2\in\{-1,0, 1\},$ denote
 \be\label{amss1}\ba  \tilde{H}(a_1,a_2)=a_1H^1+a_2H^2,\quad\tilde{u}(a_1,a_2)=a_1u^1+a_2u^2.\ea\ee
It thus follows from \eqref{a1}$_3$ that
 \be\label{amss2}\ba  \tilde{H}_t-\nu\Delta \tilde{H}=H\cdot\na\tilde{u}-u\cdot \na\tilde{H}+\tilde{H}\div u .\ea\ee
 Integrating \eqref{amss2} multiplied by $4\nu^{-1}\tilde{H} \triangle |\tilde{H}|^2$  over $\mr^2$  leads to
\be\label{wt3.51}\ba & \nu^{-1}\left(  \|\na |\tilde{H}|^2\|^2_{L^2}\right)_t+ {2}\|\Delta |\tilde{H}|^2\|^2_{L^2}\\
&=4\int  |\na \tilde{H}|^2 \Delta |\tilde{H}|^2dx
-4\nu^{-1}\int H\cdot \na \tilde{u}\cdot\tilde{H}\Delta |\tilde{H}|^2dx\\
&\quad+4\nu^{-1}\int \div u |\tilde{H}|^2 \Delta |\tilde{H}|^2dx+ 2\nu^{-1}\int u\cdot \na |\tilde{H}|^2 \Delta |\tilde{H}|^2dx \\
&\le  C  \|\na u\|^4_{L^4}+C \|\na H\|^4_{L^4}+C  \||H|^2\|^4_{L^4}
+ \|\Delta |\tilde{H}|^2\|^2_{L^2},
\ea\ee
where we have used the following simple fact that
\bnn\ba
 2\int u\cdot \na |\tilde{H}|^2 \Delta |\tilde{H}|^2dx &=
 -2\int \partial_i u\cdot \na |\tilde{H}|^2 \partial_i |\tilde{H}|^2dx
 +\int \div  u  |\na |\tilde{H}|^2|^2  dx\\
&   \leq C\|\na u\|^4_{L^4}+C\|\na H\|^4_{L^4}+C\||H|^2\|^4_{L^4}.
\ea\enn
Integrating \eqref{wt3.51} over $(0,T)$, and using  \eqref{h1}
and \eqref{es-1'}, we obtain
 \be\label{nwt4'}
 \sup_{0\le t\le T} \|\na |\tilde{H}|^2\|^2_{L^2}
 +\int_0^T \|\Delta |\tilde{H}|^2\|^2_{L^2}dt\leq C.
\ee

Noticing that
\be\label{wt3.57'1}\ba
 \||\na H||H| \|^2_{L^2}\le &  \|\na|\tilde{H}(1,0)|^2\|^2_{L^2}+ \|\na|\tilde{H}(0,1)|^2\|^2_{L^2} \\
 &+ \|\na|\tilde{H}(1,1)|^2\|^2_{L^2}+ \|\na|\tilde{H}(1,-1)|^2 \|^2_{L^2},\ea\ee
and that
\be\label{wt3.571}\ba
\||\Delta H||H| \|^2_{L^2}\le  &  C \|\na H\|^4_{L^4}+
\|\Delta |\tilde{H}(1,0)|^2\|^2_{L^2}+\|\Delta|\tilde{H}(0,1)|^2\|^2_{L^2}\\
& + \|\Delta |\tilde{H}(1,1)|^2\|^2_{L^2}
+ \|\Delta |\tilde{H}(1,-1)|^2\|^2_{L^2},
\ea\ee
then we can get \eqref{es-2} from \eqref{nwt4'}-\eqref{wt3.571}.
\hfill $\Box$
\vspace{2mm}

In order to improve the regularity estimates on $\n$, $u$ and $H$,
we start some basic energy estimates on the material derivatives of $u$.

\begin{lemma}\label{le3}
Under the condition \eqref{non-blowup}, it holds that for $0\leq T\leq T^*$,
\be\label{es-3}
\sup_{0\le t\leq T}\left(\|\n^{1/2}\dot{u}\|^2_{L^2}+
\|\nabla u\|^4_{L^4}\right)
+\int_0^T \|\nabla\dot{u}\|^2_{L^2} dt\leq C.
\ee
\end{lemma}

{\bf Proof}: We will follow an idea due to Hoff \cite{Hoff-95}. Operating $\pa/\pa_t+\div(u\cdot~)$ to $ (\ref{a1})_2^j $ and multiplying the resulting equation by $\dot{u}^j$, one gets by some simple calculations that
\be\label{wt3.40}\ba
\frac{1}{2}\left(\int \n |\dot{u}^j|^2dx\right)_t
&=\mu\int\dot{u}^j(\Delta u^j_t+\div(u\Delta u^j))dx\\
&\quad+(\mu+\lambda)\int\dot{u}^j(\pa_t\pa_j(\div u)
+\div(u\pa_j(\div u)))dx\\
&~~~-\int\dot{u}^j(\pa_jP_t+\div(u\pa_jP))dx\\
&\quad-\frac{1}{2}\int\dot{u}^j(\pa_t\pa_j|H|^2+\div(u\pa_j|H|^2))dx\\
&~~~+\int\dot{u}^j(\pa_t(H\cdot\na H^j)+\div(u(H\cdot\na H^j)))dx
\triangleq\sum^5_{i=1}J_i.
\ea\ee
First, integration by parts gives
\be\label{j1}
\ba
J_1&=\mu\int\dot{u}^j(\Delta u^j_t+\div(u\Delta u^j))dx
=-\mu\int(\partial_i\dot{u}^j\partial_i u^j_t
+\Delta u^j u\cdot\nabla \dot{u}^j) dx\\
&=-\mu\int(|\nabla\dot{u}|^2-\partial_i\dot{u}^j u^k\partial_k\partial_i u^j
-\partial_i\dot{u}^j\partial_i u^k\partial_k u^j
+\Delta u^j u\cdot\nabla \dot{u}^j) dx\\
&=-\mu\int(|\nabla\dot{u}|^2+\partial_i\dot{u}^j\partial_k u^k\partial_i u^j
-\partial_i\dot{u}^j\partial_i u^k\partial_k u^j
-\partial_i u^j\partial_i u^k\partial_k \dot{u}^j)dx\\
&\le -\frac{3\mu}{4}\|\nabla\dot{u}\|^2_{L^2}+C\|\nabla u\|^4_{L^4}.
\ea
\ee
Similarly,
\be\label{j2}
J_2\leq -\frac{\mu+\lambda}{2}\|\div\dot{u}\|^2_{L^2}
+C\|\na u\|^4_{L^4}.
\ee
It follows from integration by parts, and \eqref{p}, \eqref{es-1} that
\be\label{j3}
\ba
J_3&=-\int\dot{u}^j(\pa_jP_t+\div(u\pa_jP))dx
=\int(\pa_j\dot{u}^jP_t+\pa_jP u\cdot\nabla \dot{u}^j)dx\\
&=-\int((\gamma-1)\pa_j\dot{u}^jP\div u+\pa_j\dot{u}^j\div(Pu)
+P\pa_j(u\cdot\nabla\dot{u}^j))dx\\
&=-\int((\gamma-1)\pa_j\dot{u}^jP\div u+\pa_j\dot{u}^j\div(Pu)
+Pu\cdot\nabla\pa_j\dot{u}^j+P\pa_ju\cdot\na\dot{u}^j)dx\\
&=-\int((\gamma-1)\pa_j\dot{u}^jP\div u+P\pa_ju\cdot\nabla\dot{u}^j)dx\\
&\leq \frac{\mu}{4}\|\na \dot{u}\|^2_{L^2}+C.
\ea
\ee

Next, it follows from  \eqref{a1}$_3$ and \eqref{g1} that
\be\label{j4}\ba J_4&= \int\pa_j\dot{u}^j H\cdot H_tdx
+\frac{1}{2}\int u\cdot\na\dot{u}^j\pa_j|H|^2 dx\\
&= \frac{1}{2}\int\pa_j\dot{u}^j \div u|H|^2dx
-\frac{1}{2}\int \pa_ju\cdot\na \dot{u}^j|H|^2 dx\\
&\quad+ \int\pa_j\dot{u}^j H\cdot (H\cdot \na u+\nu \Delta H-H\div u)dx\\
&\leq C\int|\na \dot{u}| |\na u||H|^2dx
+C\int|\na \dot{u}| |\Delta H\cdot H| dx\\
&\leq \frac{\mu}{8}\|\na \dot{u}\|^2_{L^2}+C\|\na u\|^4_{L^4}
+C\|H\|^8_{L^8} +C\|\Delta H\cdot H\|^2_{L^2}.
\ea\ee
Similar to \eqref{j4}, we estimate $J_5$ as follows
\be\label{j5}\ba J_5 \leq
\frac{\mu}{8}\|\na \dot{u}\|^2_{L^2}+C\|\na u\|^4_{L^4}
+C\|H\|^8_{L^8} +C\|\Delta H\cdot H\|^2_{L^2}.
\ea\ee
Putting \eqref{j1}-\eqref{j5} into \eqref{wt3.40}, which together with \eqref{h1} gives
\be\label{wt3.46}\ba
& \frac{\rm d}{\rm dt}\|\n^{1/2}\dot{u}\|^2_{L^2}+\|\na\dot{u}\|^2_{L^2}dx
 \leq C \|\na u\|^4_{L^4}+C \||\Delta H| |H|\|^2_{L^2}+C.
\ea\ee
Integrating \eqref{wt3.46} over $(0,T)$,
using \eqref{es-1'} and \eqref{es-2}, one has
\be\label{wt-2}
\sup_{0\leq t\leq T}\|\n^{1/2}\dot{u}\|^2_{L^2}
+\int_0^T\|\na \dot{u}\|^2_{L^2}dt\leq C.
\ee
Then \eqref{es-3} can be obtained directly from \eqref{h18},
\eqref{non-blowup}, \eqref{basic}, \eqref{es-2} and \eqref{wt-2}.

\hfill $\Box$
\vspace{2mm}

Next, the following Lemma \ref{le4} combined with Lemma \ref{lemma2.6} will be
useful to estimate the $L^p$-norm of $\n\dot{u}$ and obtain the regularity
estimates on $\n$.

\begin{lemma}\label{le4}
Under the condition \eqref{non-blowup}, then there
exists a positive constant $N_1$ depending only on $N_0$, $T$ and
the initial data such that  for $0\leq T\leq T^*$,
\be\label{rho-l1}
\int_{B_{N_1}}\n(x,t)dx\geq \frac14.
\ee
\end{lemma}

{\bf Proof}: First, multiplying $\eqref{a1}_1$ by $(1+|x|^2)^{1/2}$ and
integrating the resulting equality over $\O$, we obtain after integration
by parts and using both \eqref{basic} and the fact that
\be\label{rho-1}
\int\n~ dx=\int \n_0~dx=1.
\ee
This  leads to
\be\label{wt-3}
\sup_{0\leq t\leq T}\int \n(1+|x|^2)^{1/2}dx\leq C.
\ee

Next, for $N>1$, let $\varphi_N(x)$ be a $C_0^{\infty}(\O)$ function such that
\begin{equation*}
0\leq\varphi_N(x)\leq 1,\quad
\varphi_N(x)=\left\{\begin{array}{lll}
0\quad {\rm if}~|x|\leq N,\\
1\quad {\rm if}~|x|\geq 2N,
\end{array}\right.
\quad |\nabla\varphi_N|\leq 2N^{-1}.
\end{equation*}
It follows from $\eqref{a1}_1$, \eqref{basic} and \eqref{rho-1} that
\bnn
\ba
\frac{\rm d}{\rm dt}\int\n\varphi_Ndx&=\int \n u\cdot\nabla\varphi_N dx
\geq-2N^{-1}\left(\int\n dx\right)^{1/2}\left(\int\n|u|^2 dx\right)^{1/2}
\geq-CN^{-1}
\ea
\enn
which implies
\be\label{wt-1}
\int\n\varphi_Ndx\geq\int\n_0\varphi_N dx-CN^{-1}T.
\ee
It thus follows from \eqref{oy3.8} and \eqref{wt-1} that for
$N_1\triangleq2(2+N_0+4CT)$,
$$
\int_{B_{N_1}}\n dx\geq\int\n\varphi_{\frac{N_1}{2}}dx\geq \frac14,
$$
which gives \eqref{rho-l1}. The proof of Lemma \ref{le4} is completed.
\hfill $\Box$
\vspace{2mm}

\

The next key lemma is used to bound the density gradient and
$L^1(0,T;L^{\infty})$-norm of $\na u$.

\begin{lemma}\label{le-5}
Under the condition \eqref{non-blowup}, for any $q>2$, it holds that
\be\label{differ-rho}
\ba
  &\sup_{0\le t\le T}\left(\|\n\|_{H^1\cap W^{1,q}}+\|\na u\|_{H^1}\right)
  +\int_0^T\|\na^2 u\|^2_{L^q}dt\leq C.
\ea
\ee
\end{lemma}

{\bf Proof}:  In fact, for $p\in [2,q]$, $|\na \n|^p$ satisfies
\bnnn \ba
& (|\nabla\n|^p)_t + \text{div}(|\nabla\n|^pu)+ (p-1)|\nabla\n|^p\text{div}u  \\
 &+ p|\nabla\n|^{p-2}(\nabla\n)^t \nabla u (\nabla\n) +
p\n|\nabla\n|^{p-2}\nabla\n\cdot\nabla\text{div}u = 0.\ea
\ennn
Thus, \be\la{L11}\ba
\frac{d}{dt} \norm[L^p]{\nabla\n}  &\le
 C(1+\norm[L^{\infty}]{\nabla u} )
\norm[L^p]{\nabla\n} +C\|\na^2u\|_{L^p}\\ &\le
 C(1+\norm[L^{\infty}]{\nabla u} )
\norm[L^p]{\nabla\n} +C\|\n\dot u\|_{L^p}+C\||H||\na H|\|_{L^p}, \ea\ee
due to
\be
\la{ua1}\|\na^2 u\|_{L^p}\le   C\left(\|\n\dot u\|_{L^p}+ \|\nabla
P \|_{L^p}+\||H||\na H|\|_{L^p}\right),\ee
which follows from the standard
$L^p$-estimate for the following elliptic system:
 \bnn  \mu\Delta
u+(\mu+\lambda)\na {\rm div}u=\n \dot u+\na P+\frac{1}{2}\na|H|^2-H\cdot\na H,\quad \, u\rightarrow
0\,\,\mbox{ as } |x|\rightarrow \infty. \enn

Next, the Gagliardo-Nirenberg inequality,  \eqref{evf} and
\eqref{h19}  implies
 \be\la{419}\ba  &\|\div u\|_{L^\infty}+\|\o\|_{L^\infty}  \\
 &\le C \|F\|_{L^\infty}+C\|P\|_{L^\infty}+C\|\o\|_{L^\infty}\\
 &\le C(q) +C(q) \|\na F\|_{L^q}^{q/(2(q-1))}
 +C(q) \|\na \o\|_{L^q}^{q/(2(q-1))}\\
 &\le C(q) +C(q) \left(\|\n\dot u\|_{L^q}
 +\||H||\na H|\|_{L^q}\right)^{q/(2(q-1))} ,
 \ea\ee
which, together with Lemma \ref{le9}, \eqref{ua1} and \eqref{es-1},
yields that
\be\label{b24}\ba
\|\na u\|_{L^\infty }\le& C\left(\|{\rm div}u\|_{L^\infty }+
\|\o\|_{L^\infty } \right)\log(e+\|\na^2 u\|_{L^q})
+C\|\na u\|_{L^2} +C \\\le& C\left(1+\|\n\dot u\|_{L^q}^{q/(2(q-1))}
+\||H||\na H|\|_{L^q}^{q/(2(q-1))}\right)\\
&\cdot\log(e+\|\rho \dot u\|_{L^q} +\||H||\na H| \|_{L^q}
+\|\na \rho\|_{L^q})\\
\le& C\left(1+\|\n\dot u\|_{L^q} +\||H||\na H|\|_{L^q}\right)
\log(e+   \|\na \rho\|_{L^q}) .
\ea\ee

It follows from Lemma \ref{lemma2.6}, Lemma \ref{le4},
 \eqref{wt-3} and \eqref{es-3} that,
\be\label{rho-u}
\int_0^T\|\n\dot{u}\|^2_{L^q}dt\leq C\int_0^T(\|\n^{1/2}\dot{u}\|^2_{L^2}+\|\na\dot{u}\|^2_{L^2})
\leq C.
\ee

Moreover, we have by H\"older's inequality, Gagliardo-Nirenberg
inequality that,
\be\label{wt-4}
\ba
\||H||\na H|\|_{L^q}&\leq C\|H\|_{L^{2q}}\|\na H\|_{L^{2q}}
\leq C\|H\|^{1/q}_{L^2}\|\na H\|_{L^2}\|\na^2H\|^{(q-1)/q}_{L^2}\\
&\leq C\|\na H\|_{L^2}(\|H\|_{L^2}+\|\na^2 H\|_{L^2}),
\ea\ee
integrating this over $(0,T)$, which together with \eqref{basic} and \eqref{es-1}, yields
\be\label{wt-5}
\int_0^T \||H||\na H|\|_{L^q}^2dt\leq C.
\ee

Then, substituting \eqref{b24} into \eqref{L11} where $p=q$,
we deduce from Gronwall's inequality, \eqref{rho-u} and \eqref{wt-5} that
\be \label{rho-q}
 \sup\limits_{0\le t\le T}\|\nabla \rho\|_{L^q}\le C,
\ee
which, along with \eqref{b24},  and \eqref{ua1}, shows
\be\label{mhd7}
\int_0^T\left(\|\na u\|_{L^{\infty}}+\|\na^2u\|^2_{L^q}\right)dt\leq C.
\ee

Finally, taking $p=2$ in (\ref{L11}), one gets by using \eqref{h1},
\eqref{es-1}, \eqref{mhd7}, and Gronwall's inequality that
\be\label{rho-2}
\sup\limits_{0\le t\le T}\|\nabla \n\|_{L^2}\le C.
\ee
The standard $L^2$-estimate for the elliptic system,
\eqref{es-2}, \eqref{es-3} and \eqref{rho-2} lead to
\be
\sup_{0\leq t\leq T}\|\na^2u\|^2_{L^2}\leq
C\sup_{0\leq t\leq T}\left(\|\n\dot u\|^2_{L^2}
+\|\na P\|^2_{L^2}+\||H||\na H|\|^2_{L^2}\right)\leq C,
\ee
which together with \eqref{rho-1}, \eqref{non-blowup}, \eqref{es-1}
and \eqref{rho-q}-\eqref{rho-2} finishes the proof of Lemma \ref{le-5}.
\hfill $\Box$
\vspace{2mm}

Next, it follows from  \eqref{3h}, \eqref{rho-l1}, and the Poincar\'e-type
inequality  \cite[Lemma 3.2]{F2} that for  $s>2$, $\eta\in(0,1], $
and $t\in [0,T],$
\be \label{ljo1}
\|u\bar x^{-1}\|_{L^2}+ \|u\bar x^{-\eta}\|_{L^{s/\eta}}\le
C(s,\eta)\|\n^{1/2}u\|_{L^2}+C(s,\eta)\|\na u\|_{L^2},
\ee
which together with \eqref{basic} and \eqref{es-1} gives
\be\label{u-wei'}
\|u\bar x^{-1}\|_{L^2}+\|u\bar{x}^{-\eta}\|_{L^{s/\eta}}\leq C(\eta,s).
\ee
With the help of \eqref{u-wei'}, we can get the following spatial weighted mean estimate of the density, which has been proved in \cite[Lemma 4.2]{lx1}.
\begin{lemma} \label{le6}
Under the condition \eqref{non-blowup}, it holds that for $a>1$,
$q>2$ and $0\leq T\leq T^*$,
\be \label{r-wei} \ba
&\sup_{0\le t\le T} \|  \bar x^a\n \|_{L^1\cap H^1\cap W^{1,q}}  \le C.
\ea
\ee
\end{lemma}

\begin{lemma}\label{newle}
Under the condition \eqref{non-blowup}, it holds that for $a>1$ and
$0\leq T\leq T^*$,
\be
\ba\label{gj10}
\sup_{0\le t\le T} \|H\bar{x}^{a/2}\|_{L^2}^2 +\int_{0}^{T} \|\na H\bar{x}^{a/2}\|_{L^2}^2dt\le C,
\ea
\ee
\be
\ba\label{gj10'}
\sup_{0\le t\le T} \|\na H\bar{x}^{a/2}\|_{L^2}^2
+\int_{0}^{T} \|\na^2H\bar{x}^{a/2}\|_{L^2}^2dt\le C.
\ea
\ee
\end{lemma}

{\bf Proof}:  First, multiplying \eqref{a1}$_3$  by $H\bar{x}^a$
and integrating by parts yields
\be \label{wt4.1} \ba
&\frac{1}{2}\left(\int |H|^2\bar{x}^adx\right)_t
+\nu \int |\na H|^2\bar{x}^adx=\frac{\nu}{2}\int |H|^2\Delta\bar{x}^adx\\
&\quad+\int H\cdot\na u\cdot H\bar{x}^adx
-\frac{1}{2}\int \div u|H|^2\bar{x}^adx
+\frac{1}{2}\int |H|^2u\cdot\na\bar{x}^adx
\triangleq\sum^4_{i=1}K_i.
\ea\ee

Direct calculations and \eqref{es-1} yield that
\be \label{k1}
 |K_1| \leq C\int |H|^2 \bar{x}^a \bar{x}^{-2}\log^{4}(e+|x|^2) dx\leq C\int |H|^2 \bar{x}^a dx,\ee
and that
\be\label{k2-3}\ba
 |K_2|+|K_3|&\leq C\int |\na u||H|^2 \bar{x}^a dx\leq C \|\na u\|_{L^2}\| H   \bar{x}^{a/2}\|_{L^4}^2\\
 &\leq C \|\na u\|_{L^2}\| H   \bar{x}^{a/2}\|_{L^2} (\| \na H \bar{x}^{a/2}\|_{L^2}+\| H \na\bar{x}^{a/2}\|_{L^2})\\
 &\leq C\| H \bar{x}^{a/2}\|_{L^2}^2 +\frac\nu4\| \na H \bar{x}^{a/2}\|_{L^2}^2.
\ea\ee

Then, it follows from H\"older's inequality, \eqref{g1}
and \eqref{u-wei'} that
\be\label{k4}\ba
|K_4|&\leq    C\| H \bar{x}^{a/2}\|_{L^4}\| H \bar{x}^{a/2}\|_{L^2} \|u\bar{x}^{-3/4}\|_{L^{4}}\\
&\leq   C\| H \bar{x}^{a/2}\|_{L^4}^2+C\|H \bar{x}^{a/2}\|_{L^2}^2 \|u\bar{x}^{-3/4}\|_{L^{4}}^2 \\
&\leq C\|H\bar{x}^{a/2}\|_{L^2}^2+\frac\nu4\| \na H \bar{x}^{a/2}\|_{L^2}^2.
\ea\ee
Putting \eqref{k1}-\eqref{k4} into \eqref{wt4.1},
after using Gronwall's inequality, we have
\be
\ba\label{wtnew-gj10}
\sup_{0\le t\le T} \int \bar{x}^{a}|H|^2 dx
+\int_{0}^{T} \int  \bar{x}^{a}|\na H|^2 dxdt\le C.
\ea
\ee

Next, multiplying  \eqref{a1}$_3$ by $\Delta H\bar{x}^a$, integrating the
resultant equation by parts over $\mr^2$, it follows from
the similar arguments as \eqref{mm1} that
\be\la{AMSS5}\ba
&\frac{1}{2}\left(\int |\na H|^2\bar{x}^adx\right)_t+\nu \int |\Delta H|^2\bar{x}^adx\\
\le& C\int|\na H| |H| |\na u| |\na\bar{x}^a|dx+C\int|\na H|^2|u| |\na\bar{x}^a|dx
+C\int|\na H| |\Delta H| |\na\bar{x}^a|dx\\
&+C\int |H||\na u||\Delta H|\bar{x}^adx+C\int |\na u||\na H|^2  \bar{x}^adx
\triangleq \sum_{i=1}^5 L_i.
\ea\ee
Using Gagliardo-Nirenberg inequality, \eqref{es-3} and \eqref{gj10},
it holds that
\be\label{l1}\ba
L_1\le &C\int|\na H| |H| |\na u| \bar{x}^{a} (\bar{x}^{-1} |\na\bar{x}|)dx\\
\le & C\|H\bar{x}^{a/2}\|_{L^4}^4+C\|\na u\|_{L^4}^4
+C\|\na H\bar{x}^{a/2}\|_{L^2}^2\\
\le & C\|H\bar{x}^{a/2}\|_{L^2}^2\left(\|\na H\bar{x}^{a/2}\|_{L^2}^2+\|H\bar{x}^{a/2}\|_{L^2}^2\right)+C
+C\|\na H\bar{x}^{a/2}\|_{L^2}^2\\
\le &C+C\|\na H\bar{x}^{a/2}\|_{L^2}^2,
\ea\ee
\be\label{l2}\ba
L_2&\leq C\int |\na H|^{(4a-1)/(2a)}\bar{x}^{(4a-1)/4}
|\na H|^{1/(2a)}|u|\bar{x}^{-1/2} \bar{x}^{-1/4} |\na\bar{x}|dx\\
&\leq  C\||\na H|^{(4a-1)/(2a)}\bar{x}^{(4a-1)/4}\|_{L^{\frac{4a}{4a-1}}}
\||\na H|^{1/(2a)}\|_{L^{8a}}\|u\bar{x}^{-1/2}\|_{L^{8a}}\\
&\leq  C\|\na H \bar{x}^{a/2} \|_{L^2}^2+  C\|\na H \|_{L^4}^2\\
&\leq  C\|\na H \bar{x}^{a/2}\|_{L^2}^2+C\|\na^2 H\|_{L^2}^2+C,
\ea\ee
\be\label{l3-4}\ba
L_3+L_4 &\le \|\Delta H\bar{x}^{a/2}\|_{L^2}\|\na H\bar{x}^{a/2}\|_{L^2}
+\|\Delta H\bar{x}^{a/2}\|_{L^2}\| H\bar{x}^{a/2}\|_{L^4}\|\na u\|_{L^4}\\
&\leq \frac{\nu}{4}\|\Delta H\bar{x}^{a/2}\|^2_{L^2}
+C\|\na H \bar{x}^{a/2}\|_{L^2}^2+C,
\ea\ee
\be\label{l5}\ba
L_5\le &C\|\na u\|_{L^\infty} \|\na H\bar{x}^{a/2}\|_{L^2}^2\\
\le &C \|\na u\|_{L^2}^{(q-2)/(2q-2)}\|\na^2 u\|_{L^q}^{q/(2q-2)} \|\na H\bar{x}^{a/2}\|_{L^2}^2\\
\le &C(1+\|\na^2 u\|_{L^q}^{2})\|\na H\bar{x}^{a/2}\|_{L^2}^2.
\ea\ee
Noticing the fact that
\be\label{dj}
\ba
\int|\na^2H|^2\bar{x}^adx&=\int|\Delta H|^2\bar{x}^adx
-\int\pa_i\pa_kH\cdot\pa_kH\pa_i\bar{x}^a dx\\
&\quad+\int\pa_i\pa_iH\cdot\pa_kH\pa_k\bar{x}^a dx\\
&\leq \int|\Delta H|^2\bar{x}^adx+\frac{1}{2}\int|\na^2H|^2\bar{x}^adx
+C\int|\na H|^2\bar{x}^a dx.
\ea
\ee
Submitting \eqref{l1}-\eqref{l5} into \eqref{AMSS5}, and using \eqref{dj},
we obtain
\be\label{AMSS10}\ba
 &\left(\int |\na H|^2\bar{x}^adx\right)_t
 +\int |\na^2 H|^2\bar{x}^adx\\
&\quad\le  C(1+\|\na^2 u\|_{L^q}^{2})\|\na H\bar{x}^{a/2}\|_{L^2}^2+C(\|\na^2 H\|_{L^2}^2+1),
\ea\ee
which together with  Gronwall's inequality,  \eqref{es-1}
and \eqref{differ-rho} yields  \eqref{gj10'}.
The proof of Lemma \ref{newle} is finished.
\hfill $\Box$
\vspace{2mm}

\begin{lemma}\label{le7}
Under the condition $\eqref{non-blowup}$, it holds that for
$0\le T\le T^*$,
\begin{equation}\label{gj13}
\ba
 &\sup_{0\leq t\leq T }\left(\|\n^{1/2}u_t\|^2_{L^2}+\|H_t\|^2_{L^2}
 +\|\na H\|^2_{H^1} \right)\\
 &+\int_0^T\left(\|\na u_t\|_{L^2}^2+\|\na H_t\|_{L^2}^2
 +\|\na^2H\|^2_{L^q}\right)dt\le C.
\ea
\end{equation}
\end{lemma}

{\bf Proof}:
First, the combination of \eqref{r-wei} with \eqref{u-wei'} gives  that
for any $\eta\in(0,1]$ and any $s>2$,
\be\label{5.d2}
 \|\n^\eta u \|_{L^{s/\eta}}+ \|u\bar x^{-\eta}\|_{L^{s/\eta}}\le C(\eta,s).
\ee

Differentiating $\eqref{a1}_2$ with respect to $t$ gives
 \be\label{zb1}\ba
 &\n u_{tt}+\n u\cdot \na u_t-\mu\Delta u_t-( \mu+\lm)\na  \div u_t  \\ &=-\n_t(u_t+u\cdot\na u)-\n u_t\cdot\na u -\na P_t
 +\left(H\cdot\na H-\frac{1}{2}\na |H|^2\right)_t.
\ea\ee
Multiplying \eqref{zb1} by $u_t$, then integrating  over $\O,$
we obtain after using  $\eqref{a1}_1$ that
\be\ba  \label{na8}&\frac{1}{2}\frac{\rm d}{\rm dt} \int \n |u_t|^2dx
+\int \left(\mu|\na u_t|^2+( \mu+\lm)(\div u_t)^2  \right)dx\\
&=-2\int \n u \cdot \na  u_t\cdot u_tdx
-\int \n u \cdot\na (u\cdot\na u\cdot u_t)dx\\
&\quad-\int \n u_t \cdot\na u \cdot  u_tdx+\int P_{t}\div u_{t} dx
+\int\left(H\cdot\na H-\frac{1}{2}\na |H|^2\right)_t u_tdx\\
 &\triangleq \sum_{i=1}^{5}\bar{J}_i.
  \ea\ee
Similar to the proof of  \cite[Lemma 4.3]{lx1},
and using \eqref{differ-rho}, we have
 \be\label{ji'}\ba
 \sum_{i=1}^{4}\bar{J}_i &\le  \frac{\mu}{4}\| \na u_{t}\|_{L^{2}}^{2}
 +C \left(\| \na^{2} u \|_{L^{2}}^{2} +
 \|\n^{1/2}u_{t}\|_{L^{2}}^{2}+1\right)\\
 &\leq\frac{\mu}{4}\| \na u_{t}\|_{L^{2}}^{2}
 +C \left(\|\n^{1/2}u_{t}\|_{L^{2}}^{2}+1\right).
 \ea\ee
For the term $\bar{J}_5$, we obtain after integration by parts that
\be\label{j5'}
\ba
\bar{J}_5&=-\int H_t\cdot\na u_t\cdot H dx-\int H\cdot\na u_t\cdot H_tdx
+\int H\cdot H_t \div u_tdx\\
&\leq C\|H\|_{L^4}\|H_t\|_{L^4}\|\na u_t\|_{L^2}\\
&\le \frac{\mu}{4}\|\na u_t\|^2_{L^2}+\de\|\na H_t\|^2_{L^2}
+C(\de)\|H_t\|^2_{L^2}.
\ea\ee
Substituting \eqref{ji'} and \eqref{j5'} into \eqref{na8} leads to
\be\label{r-ut}
\ba
\frac{\rm d}{\rm dt}\|\n^{1/2}u_t\|^2_{L^2}+\|\na u_t\|^2_{L^2}
\le C\de\|\na H_t\|^2_{L^2}+C(\de)\left(\|\n^{1/2}u_t\|^2_{L^2}
+\|H_t\|^2_{L^2}+1\right).
\ea\ee

Next, differentiating $\eqref{a1}_3$ with respect to $t$ shows
  \be\la{wt4.12}\ba
H_{tt}-H_t\cdot\na u-H\cdot\na u_t+u_t\cdot\na H+u\cdot\na H_t+H_t\div u+H\div u_t=\nu \Delta H_t.
\ea\ee
 Multiplying \eqref{wt4.12} by $H_t$,  then integrating the resulting equation over $\O,$
 and using \eqref{h1}, \eqref{es-1}, \eqref{ljo1} and \eqref{gj10}, we have
\be\la{wt4.13}\ba
&\frac{1}{2}\frac{\rm d}{\rm dt} \int |H_t|^2dx+\nu \int  |\na H_t|^2dx
=\int H_t\cdot\na u\cdot H_tdx\\&\quad-\frac12\int \div u|H_t|^2dx
+\int H \cdot \na u_t\cdot H_tdx+\int u_t\cdot\na H_t\cdot Hdx\\
&\le C\|H_t\|^2_{L^4}\|\na u\|_{L^2}
+C\|H\|_{L^4}\|H_t\|_{L^4}\|\na u_t\|_{L^2}\\
&\quad+C\||H|^{\frac{1}{2a}}\|_{L^{8a}}
\|(H\bar{x}^{\frac{a}{2}})^{\frac{2a-1}{2a}}\|_{L^{\frac{4a}{2a-1}}}
\|u_t\bar{x}^{-\frac{2a-1}{4}}\|_{L^{8a}}\|\na H_t\|_{L^2}\\
& \leq \de\|\na H_t\|^2_{L^2}+C(\de)\left(\|H_t\|^2_{L^2}
+\|\n^{1/2}u_t\|^2_{L^2}+\|\na u_t\|^2_{L^2}\right),
\ea
\ee
which gives
\be\label{h-t}
\frac{\rm d}{\rm dt}\|H_t\|^2_{L^2}+\|\na H_t\|^2_{L^2}
\le C\|H_t\|^2_{L^2}+C\|\n^{1/2}u_t\|^2_{L^2}+C_2\|\na u_t\|^2_{L^2}.
\ee
Then adding \eqref{r-ut} multiplied by $C_2+2$ to \eqref{h-t}, choosing
$\de$ suitable small, we obtain
\be\label{ut-ht}
\ba
\frac{\rm d}{\rm dt}\left((C_2+2)\|\n^{1/2}u_t\|^2_{L^2}
+\|H_t\|^2_{L^2}\right)+\|\na u_t\|^2_{L^2}+\frac{1}{2}\|\na H_t\|^2_{L^2}\\
\le C\left(\|\n^{1/2}u_t\|^2_{L^2}+\|H_t\|^2_{L^2}+1\right),
\ea
\ee
which together with Gronwall's inequality  yields
\be\label{ut-ht'}
 \sup_{0\le t\le T}\left(\|\n^{1/2}u_t\|^2_{L^2}+\|H_t\|^2_{L^2}\right)
+\int_0^T\left(\|\na u_t\|^2_{L^2}+\|\na H_t\|^2_{L^2}\right)dt\leq C.
\ee

Finally, it follows from \eqref{ut-ht'}, \eqref{h1}, \eqref{u-wei'},
\eqref{es-3}, \eqref{gj10'} and Gagliardo-Nirenberg inequality that
\be\ba\label{dao-H}\|\na^2 H\|^2_{L^2}
&\leq C\||u| |\na H|\|_{L^2}^2+C\||H| |\na u|\|^2_{L^2}+C\|H_t\|^2_{L^2}\\
&\le C\|u\bar{x}^{-a/4}\|^2_{L^8}\||\na H|^{1/2}\bar{x}^{a/4}\|^2_{L^4}
\||\na H|^{1/2}\|^2_{L^8}\\
&\quad+C\|H\|^2_{L^4}\|\na u\|^2_{L^4}+C\\
&\leq C\|\na H\bar{x}^{a/2}\|_{L^2}\|\na H\|_{L^4}+C\\
&\leq \frac{1}{2}\|\na^2 H\|^2_{L^2}+C,
\ea\ee
and
\be\label{daohq}
\ba
\int_0^T\|\na^2H\|^2_{L^q}dt&\leq C\int_0^T\left(
\||u| |\na H|\|_{L^q}^2+\||H| |\na u|\|^2_{L^q}+\|H_t\|^2_{L^q}\right)dt\\
&\leq C\int_0^T\big(\|u\bar{x}^{-a/2}\|^2_{L^{2q}}\|\na H\bar{x}^{a/2}\|^2_{L^{2q}}
+\|H\|^2_{L^{2q}}\|\na u\|^2_{L^{2q}}\\
&\qquad\qquad+\|H_t\|^{4/q}_{L^{2}}\|\na H_t\|^{2(q-2)/q}_{L^{2}}\big)dt\\
&\leq C\int_0^T\big(\|\na H\bar{x}^{a/2}\|^{2/q}_{L^2}
(\|\na^2 H\bar{x}^{a/2}\|_{L^2}+\|\na H\na\bar{x}^{a/2}\|_{L^2})^{2(q-1)/q}\\
&\qquad\qquad+\|\na u\|^2_{H^1}+\|\na H_t\|^2_{L^2}+1\big) dt\\
&\leq C\int_0^T\big(\|\na^2 H\bar{x}^{a/2}\|^2_{L^2}
+\|\na H_t\|^2_{L^2}+1\big)dt\leq C,
\ea
\ee
which combined with \eqref{dao-H} leads to \eqref{gj13}.
Thus the proof of Lemma \ref{le7} is finished.
\hfill $\Box$
\vspace{2mm}

Now we are in a position to prove Theorem \ref{th1}.

{\bf Proof of Theorem \ref{th1}}: Suppose that \eqref{blowup} were false,
that is, \eqref{non-blowup} holds. Note that the generic constant $C$
in Lemma \ref{basic} and  Lemma \ref{le-5}-\ref{le7} remains uniformly
bounded for all $T<T^*$, so the functions
$(\n,u,H)\triangleq\lim_{t\rightarrow T^*}(\n,u,H)(x,t)$ satisfy the
conditions imposed on the initial data \eqref{co1} at the time $t=T^*$.
Furthermore, standard arguments yield that $\n\dot u\in C([0,T];L^2)$,
which implies
\bnn
(\n\dot u)(x,T^*)=\lim_{t\rightarrow T^*}(\n \dot{u})\in L^2.
\enn
Hence,
\bnn
-\mu\Delta u-(\mu+\lambda)\na \div u+\na P-(\na\times H)\times H|_{t=T^*}
=\sqrt{\n}(x,T^*)g(x),
\enn
with
\bnn
g(x)=\left\{
\ba
&\n(x,T^*)^{-1/2}(\n\dot u)(x,T^*),\quad &{\rm for}~~
x\in\{x|\n(x,T^*)>0\},\\
&0, \quad &{\rm for}~~x\in\{x|\n(x,T^*)=0\},
\ea
\right.
\enn
satisfying $g\in L^2$ due to \eqref{es-3}. Thus, $(\n, u)(x,T^*)$ satisfies
\eqref{cc} also. Therefore, one can take $(\n, u,H)(x, T^*)$ as
the initial data and apply Lemma \ref{th0} to extend the local
strong solution beyond $T^*$. This contradicts the assumption on $T^*$.
We thus finish the proof of Theorem \ref{th1}.





\textbf{Acknowledgments.} The author would like to express great gratitude to Prof. Feimin Huang for his valuable
comments and suggestions which greatly improved the presentation of the paper.








\begin{thebibliography}{99}



\bibitem{df2006} B. Ducomet, E. Feireisl, The equations of magnetohydrodynamics: On the interaction between matter and radiation in the
evolution of gaseous stars, \textit{Comm. Math. Phys.} \textbf{266}(2006), 595--629.



\bibitem{cw} Chen, G. Q.; Wang, D. Global solution of nonlinear magnetohydrodynamics with large initial data, \textit{J. Differential Equations,} \textbf{182} (2002), 344--376.

\bibitem{cw2} Chen, G. Q.; Wang, D.  Existence and continuous dependence of large solutions for the magnetohydrodynamic equations, \textit{Z. Angew. Math. Phys.}, \textbf{54} (2003), 608--632.

\bibitem{fjn} Fan, J.; Jiang,  S.; Nakamura, G.  Vanishing shear viscosity limit in the magnetohydrodynamic
equations, \textit{Commun. Math. Phys.}, \textbf{270}(2007), 691--708.


\bibitem{fy2} Fan, J. ; Yu, W. Strong solution to the compressible MHD equations with vacuum,
\textit{Nonlinear
Anal. Real World Appl.}, \textbf{10}(2009),  392--409.

\bibitem{F2}  Feireisl, E.  Dynamics of viscous compressible fluids. Oxford
University Press, New York,   2004.

\bibitem{F1} Feireisl, E.; Novotny, A.; Petzeltov\'{a}, H.  On the existence of globally defined weak solutions to the Navier-Stokes equations.
\textit{J. Math. Fluid Mech.}  \textbf{3} (2001),  no. 4, 358--392.




\bibitem{he-xin}
He, C.; Xin, Z. P.: On the regularity of weak solutions to the magnetohydrodynamic equations.
\textit{J. Diff.Eq.}, \textbf{213} (2005), 235--254.

\bibitem{Hoff-95}
Hoff, D. Global solutions of the Navier-Stokes equations for
multidimensional compressible flow with discontinuous initial data.
\textit{J. Differ. Eqs.}   \textbf{120} (1995), no. 1, 215--254.







\bibitem{hw1} Hu, X.; Wang, D. Global solutions to the three-dimensional full compressible magnetohydrodynamic
flows, \textit{Commun. Math. Phys.}, \textbf{283} (2008) 255--284.

\bibitem{hw2} Hu, X.; Wang, D. Global existence and large-time behavior of solutions to the threedimensional equations of compressible magnetohydrodynamic flows, \textit{Arch. Ration. Mech. Anal.}, \textbf{197} (2010) 203--238.




\bibitem{hx2} Huang, X. D.; Li, J.; Xin Z. P.
Blowup criterion for viscous barotropic flows with vacuum states.
\textit{Comm. Math. Phys.}  \textbf{301} (2011), no. 1, 23--35.

\bibitem{hlx} Huang, X. D.; Li, J.; Xin Z. P.
Serrin type criterion for the three-dimensional compressible flows.
 \textit{SIAM J. Math. Anal.},  \textbf{ 43} (2011), no. 4, 1872--1886.


\bibitem{hlx1} Huang, X. D.; Li, J.; Xin, Z. P.    Global well-posedness of classical solutions with large
oscillations and vacuum to the three-dimensional isentropic
compressible Navier-Stokes equations.   \textit{Comm. Pure Appl.
Math. } {\bf65} (2012), 549--585.

\bibitem{hlw} Huang, X. D.; Li, J.; Wang, Y. Serrin-type blowup criterion for full compressible Navier-Stokes system. \textit{Arch. Ration. Mech. Anal.}, \textbf{207}(2013), 303--316.

\bibitem{hl} Huang, X. D.; Li, J. Serrin-Type Blowup Criterion for Viscous, Compressible, and Heat Conducting Navier-Stokes and Magnetohydrodynamic Flows. \textit{Comm. Math. Phys.}, \textbf{324}(2013), 147--171.



\bibitem{ko} S. Kawashima; M. Okada, Smooth global solutions for the one-dimensional equations in magnetohydrodynamics, \textit{Proc. Japan Acad. Ser. A Math. Sci.}, \textbf{58} (1982) 384--387.

\bibitem{ka} S. Kawashima, Systems of a hyperbolic-parabolic composite type, with applications to the equations of magnetohydrodynamics,
PhD thesis, Kyoto University, 1983.

\bibitem{li-liang}
Li, J.; Liang, Z. L.
On local classical solutions to the Cauchy problem of the two-dimensional barotropic compressible Navier-Stokes equations with vacuum.
\textit{J. Math. Pures Appl.}, \textbf{102} (2014), 640--671.

\bibitem{lx1} Li, J.; Xin, Z. Global well-posedness and decay asymptotic behavior of classical solution to the compressible Navier-Stokes equations with vacuum, http://arxiv.org/abs/1310.1673.


\bibitem{lxz} Li, H.; Xu, X.; Zhang, J. Global classical solutions to 3D compressible magnetohydrodynamic equations with large oscillations and vacuum, \textit{SIAM J. Math. Anal.} (3) {\bf 45}(2013), 1356--1387.

\bibitem{L2} Lions,  P. L.  {Mathematical topics in fluid mechanics. Vol. {\bf 1}. Incompressible models.}
Oxford University Press, New York, 1996.

\bibitem{L1} Lions, P. L.    {Mathematical topics in fluid mechanics. Vol. 2. Compressible models.}  Oxford
University Press, New York,  1998.


\bibitem{chl}  Lv, B. Q.; Huang, B. On strong solutions to the cauchy problem of the 2-D compressible MHD equations with vacuum.
Accepted by Nonlinearity.

\bibitem{lv-shi-xu}
Lv, B. Q.; Shi, X. D.; Xu, X. Y. Global well-posedness and large time
asymptotic behavior of strong solutions to the 2-D compressible
magnetohydrodynamic equations with vacuum.
http://arxiv.org/abs/1402.4851.



\bibitem{nir}
 Nirenberg, L.  On elliptic partial differential equations.  \textit{Ann. Scuola Norm. Sup. Pisa}, (3){\bf 13}(1959), 115--162.





\bibitem{sun-zhang}
Sun, Y. Z; Zhang, Z. F. A blow-up criterion of strong solutions to the 2D compressible Navier-Stokes equations.
\textit{Sci. China Math.}, \text{54} (2011),  105--116.

\bibitem{swz-1}
Sun, Y.Z.; Wang, C.; Zhang, Z.F.
A Beale-Kato-Majda Blow-up criterion for the 3-D compressible
Navier-Stokes equations.
\textit{J. Math. Pures Appl.}, \textbf{95} (2011), 36--47.


\bibitem{uks} T. Umeda, S. Kawashima, and Y. Shizuta,
On the decay of solutions to the linearized equations of electromagnetofluid dynamics,
\textit{Japan J. Appl. Math.}, \textbf{1} (1984), 435--457.

\bibitem{vk} A. I. Vol'pert; S. I. Khudiaev,
On the Cauchy problem for composite systems of nonlinear equations,
\textit{Mat. Sb.}, \textbf{87}(1972), 504--528.

\bibitem{w} Wang, D.
Large solutions to the initial-boundary value problem for planar
magnetohydrodynamics,
\textit{SIAM J. Appl. Math.}, \textbf{63} (2003), 1424--1441.

%

\bibitem{xu-zhang}
Xu, X. Y.; Zhang, J. W. A blow-up criterion for 3D compressible
 magnetohydrodynamic equations with vaccum.
\textit{Math. Models Method. Appl. Sci.},  \textbf{22} (2012), 1150010.


 \end{thebibliography}
\end{document}